\documentclass[journal, letter, 9pt]{IEEEtran}
\usepackage{amsmath}
\usepackage{amsfonts}
\usepackage{graphics}
\usepackage{pstricks}
\usepackage{array}
\usepackage{epsf}
\usepackage{graphicx}
\usepackage{multirow}
\usepackage{cite}
\usepackage{subcaption}
\usepackage[font=small]{caption}

\newcommand{\bfg}[1]{\boldsymbol{#1}}
\newcommand{\bfp}[1]{\boldsymbol{#1}'}
\newcommand{\bfpp}[1]{\boldsymbol{#1}''}
\newcommand{\bfd}[1]{\dot{\boldsymbol{#1}}}
\newcommand{\bfdd}[1]{\ddot{\boldsymbol{#1}}}
\newcommand{\bfb}[1]{\boldsymbol{\rm #1}}
\newcommand{\flux}{\varphi}
\newcommand{\ii}{\imath}
\newcommand{\jj}{\jmath}
\newcommand{\geom}[1]{\bfg #1 \bfp #1}
\newcommand{\gw}[1]{\bfb \Omega_{#1}}
\newcommand{\gr}[1]{\frac{\geom{#1}}{#1^2}}
\newcommand{\gf}[1]{\rho_{#1} + \gw{#1}}
\newcommand{\inner}[1]{\frac{\bfg #1}{#1} \cdot \frac{\bfp #1}{#1}}
\newcommand{\out}[1]{\frac{\bfg #1}{#1} \wedge \frac{\bfp #1}{#1}}
\newcommand{\e}[1]{\boldsymbol{\rm e}_{#1}}
\newcommand{\ep}[1]{\boldsymbol{\rm e}'_{#1}}

\begin{document}

\title{A Geometrical Interpretation of Frequency}

\author{Federico~Milano,~\IEEEmembership{Fellow,~IEEE}%
\thanks{F.~Milano is with School of Electrical and Electronic
  Engineering, University College Dublin, Belfield Campus, Dublin 4,
  Ireland.  E-mails: federico.milano@ucd.ie}
\thanks{This work is supported by the European Commission by funding
  F.~Milano under project EdgeFLEX, Grant No.~883710.}}

\maketitle

\begin{abstract}
  The letter provides a geometrical interpretation of frequency in
  electric circuits.  According to this interpretation, the frequency
  is defined as a multivector with symmetric and antisymmetric
  components.  The conventional definition of frequency is shown to be
  a special case of the proposed theoretical framework.  Several
  examples serve to show the features, generality as well as
  practical aspects of the proposed approach.
\end{abstract}


\IEEEpeerreviewmaketitle

\section{Introduction}
\label{sec:intro}

In power system applications, the frequency of an ac signal is
conventionally defined as the time derivative of the argument of the
cosine function of the signal itself \cite{IEEE118}.  This definition
appears to have some issues.  First, it depends on the representation
of the signal itself.  This has led to a tremendous number of
publications, each of which using as starting point a different
representation \cite{8586583}.  Second, the value of the frequency
often depends on the transformation utilized to represent the ac
signal.  A good criterion to decide if a transformation is robust is
to check whether a signal can be fully reconstructed to its original
state if the inverse transformation is applied to the transformed
signal \cite{Paolone:2020}.  This is a sensible criterion but does not
guarantee the correctness and consistency of the estimation of the
frequency itself.  The value of the estimated frequency should be
always the same (invariant) independently from the transformation.  A
third issue with the common definition of frequency and a large number
of existing techniques to estimate the frequency is that they do not
account for variations of the magnitude of the signal.  
This assumption poses serious issues for the estimation of the
frequency from measurements.  The conventional definition, in fact,
implicitly assumes that one is able to measure the phase angle
independently from the magnitude and that the measured signal is a
sine wave.  But this is not always the case, in particular, in
transient conditions.   The theoretical framework and
definition of generalized frequency proposed in this letter address
the issues above.

\section{Outlines of Vector Operations and Space Curves} 
\label{sec:geometry}

We first provide some outlines of operations with vectors and space
curves.   In the remainder of the letter, vectors are
indicated in bold face (e.g., $\bfg v$), whereas scalar quantities are
in normal face (e.g., $v$).   Let
$\bfg x = (x_1, x_2, \dots, x_n)$ and
$\bfg y = (y_1, y_2, \dots, y_n)$ be two $n$-dimensional vectors in
$\mathbb{R}^n$.

The \textit{inner product} is defined as:
\begin{equation}
  \label{eq:inner}
  \bfg x \cdot \bfg y = \sum_{i=1}^n x_i y_i \, .
\end{equation}
For example, in $\mathbb{R}^3$,
$\bfg x \cdot \bfg y = x_1y_1 + x_2y_2 + x_3y_3$. The inner product is
symmetric, associative, and commutative.  In particular, the inner
product of a vector by itself gives:
\begin{equation}
  x = |\bfg x| = \sqrt{\bfg x \cdot \bfg x} \, ,
\end{equation}
 where $x$ is the magnitude of $\bfg x$. 

The \textit{outer product} is defined as:
\begin{equation}
  \label{eq:outer}
  \bfg x \otimes \bfg y =
  \begin{bmatrix}
    x_1 y_1 & \dots & x_1 y_n \\
    \vdots & \ddots & \vdots \\
    x_n y_1 & \dots & x_n y_n     
  \end{bmatrix} . 
\end{equation}

The \textit{wedge product} is defined as: 
\begin{equation}
  \label{eq:wedge}
  \bfg x \wedge \bfg y = \bfg x \otimes \bfg y - \bfg y \otimes \bfg x \, .
\end{equation}
For example, in $\mathbb{R}^3$, the wedge product gives:
\begin{equation}
  \bfg x \wedge \bfg y  =
  \begin{bmatrix}
    0 & b_{12} & -b_{31 }\\
    -b_{12} & 0 & b_{23} \\
    b_{31} & -b_{23} & 0
  \end{bmatrix} ,
\end{equation}
where $b_{ij} = x_iy_j - y_ix_j$.  The result of the wedge product is
a \textit{bivector}.  In the remainder of this letter it will be
indicated with an uppercase bold symbol, e.g.,
$\bfb B = \bfg x \wedge \bfg y$, where $\bfb B$ is a skew-symmetric
matrix.  The wedge product is antisymmetric, associative, and
anti-commutative.  The latter means that
$\bfg x \wedge \bfg y = - \bfg y \wedge \bfg x$ and, consequently
$\bfg x \wedge \bfg x = \bfg 0$.
In $\mathbb{R}^3$, the wedge product is similar to the \textit{cross
  product} $\bfg x \times \bfg y$, although the result of the cross
product is a vector not a tensor.  For the developments of this
letter, it is relevant to note that in Euclidean metric, the magnitude
of a bivector is given by:
\begin{equation}
  |\bfg x \wedge \bfg y| = |\bfb B| = \sqrt{\sum_{i=1}^n \sum_{j>i}^n b^2_{ij}} \, .
\end{equation}
For example, in $\mathbb{R}^3$:
\begin{equation}
  |\bfg x \wedge \bfg y| = \sqrt{b_{12}^2 + b_{23}^2 + b_{31}^2} \, .
\end{equation}
%
%
The \textit{geometric product} is defined as:
\begin{equation}
  \label{eq:gp}
  \bfg x \bfg y = \bfg x \cdot \bfg y + \bfg x \wedge \bfg y \, .
\end{equation}
The result of the geometric product, which is called
\textit{multivector}, consists of two components.  The first
component, $\bfg x \cdot \bfg y$, is a scalar that represents the
projection of $\bfg y$ onto the vector $\bfg x$.  The second
component, $\bfg x \wedge \bfg y$, represents a bivector orthogonal to
the space defined by the vectors $\bfg x$ and $\bfg y$.  It may seem
strange at first to sum a scalar with a bivector but this is
exactly the same kind of operation that is intended when one writes a
complex number as $a + \jj b$.  Section \ref{sec:examples} shows that,
in fact, the algebra of complex numbers is a special case of the
algebra of multivectors.

In this work, we are interested in time-dependent $n$-dimensional
\textit{curves} (or trajectories), i.e.,
$\bfg x(t) = (x_1(t), x_2(t), \dots, x_n(t))$, where $t$ is time.  The
time derivative of $\bfg x$ is defined as:
\begin{equation}
  \bfp x  = \frac{d \bfg x}{dt}  =
  (x'_1, x'_2, \dots, x'_n) \, .
\end{equation}
From the geometrical point of view, $\bfp x$ is the \textit{tangent
  vector} of the curve $\bfg x$.  Let us define $s$ as the \textit{arc
  length} of the curve $\bfg x$, then the following property holds:
\begin{equation}
  \label{eq:length}
  s' = \frac{d s}{dt}  =
  \sqrt{\bfp x \cdot \bfp x} = x' \, .
\end{equation}
It is important to note that the arc length $s$ and thus also its
derivatives are invariant with respect to the system of coordinates.

It is relevant to define the derivative of the curve $\bfg x$ with
respect to $s$, as follows:
\begin{equation}
  \label{eq:xdot}
  \bfd x = \frac{d \bfg x}{ds} =
  \frac{d \bfg x}{dt} \frac{dt}{ds} =
  \frac{\bfp x}{x'} \, ,
\end{equation}
where we have used \eqref{eq:length} and the identity $dt/ds = 1/s'$.
From \eqref{eq:xdot}, it follows that $\bfd x \cdot \bfd x = 1$.  The
vector $\bfd x$ is tangent to $\bfg x$ and that the tangent vector to
a curve is the unit vector if the arc length is chosen as a parameter.

Finally, it is relevant to recall the definition of another invariant
quantity, namely the \textit{curvature}, that plays an important role in
differential geometry.  The curvature is defined as:
\begin{equation}
  \label{eq:kappa}
  \kappa = |\bfd x \wedge \bfdd x| = \frac{|\bfp x \wedge \bfpp x|}{(x')^{3}} \, ,
\end{equation}
where
\begin{equation}
  \label{eq:xdotdot}
  \bfdd x = \frac{d \bfd x}{ds} = \frac{\bfpp x}{(x')^2} -
  \frac{x'' \, \bfp x}{(x')^3}
\end{equation}
is the tangent vector to $\bfd x$ and satisfies the condition
$\bfd x \cdot \bfdd x = 0$.

We are now ready to present the main contribution of this work.

\section{Frequency as a Multivector}

Let us start with the vector of the magnetic flux, $\bfg \flux$.
According to the Faraday's law of induction, one has:
\begin{equation}
  \label{eq:flux}
  -\bfp \flux = \bfg v \, ,
\end{equation}
where $\bfg v$ is the vector of the voltage and
 the minus accounts for the Lenz's law but is not crucial
for the discussion below.  On the other hand, it is important to note
that $\bfg \flux$ does not need to be known or to be measurable.  In
the context of this work, $\bfg \flux$ serves only to define the
macroscopic effect of the magnetic field.  In this context, the most
important property of $\bfg \flux$ is that its time derivative is the
vector of the voltage.  If one interprets the components of the vector
of the flux as the coordinates of a \textit{curve}, say
$\bfg x = -\bfg \flux$, then the voltage $\bfg v = \bfp x$ is the
tangent vector to this curve.

According to the definitions given in Section \ref{sec:geometry}, one
has:
\begin{equation}
  \label{eq:v}
  s' = |\bfg \flux' \cdot \bfg \flux'| =
  \flux' = v \, ,
\end{equation}
and
\begin{equation}
  \label{fluxdot}
  \bfd x =
  -\bfd \flux =
  \frac{\bfp \flux}{s'} =
  \frac{\bfg v}{v} \, ,
\end{equation}
and
\begin{equation}
  \label{eq:fluxdotdot}
  \bfdd x =
  -\bfdd \flux =
  \frac{\bfp v}{v^2} -
  \frac{v' \, \bfg v}{v^3} \, .
\end{equation}
Since $\bfd x \cdot \bfdd x = 0$, one obtains:
\begin{equation}
  0 = \bfd \flux \cdot \bfdd \flux =
  \frac{\bfg v \cdot \bfp v}{v^3} -
  \frac{v' \, \bfg v \cdot \bfg v }{v^4}
\end{equation}
which leads to:
\begin{equation}
  \label{eq:rho}
  \rho_v = \frac{v'}{v} = \inner{v} \, .
\end{equation}

Similarly, from the definition of curvature in \eqref{eq:kappa}, one
obtains:
\begin{equation}
  \kappa_v = | \bfd \flux \wedge \bfdd \flux | =
  \left | \frac{\bfg v}{v} \wedge \left ( \frac{\bfp v}{v^2} -
  \frac{v' \, \bfg v }{v^3} \right )  \right | ,
\end{equation}
and remembering that $\bfg v \wedge \bfg v = 0$, one obtains:
\begin{equation}
  \label{eq:kappabis}
  \kappa_v = 
  \left | \frac{\bfg v}{v} \wedge \frac{\bfp v}{v^2} \right | ,
\end{equation}
We define the magnitude of the frequency of $\bfg v$, say
$\omega_v$, as:
\begin{equation}
  \label{eq:omega}
  \omega_v = v \, \kappa_v \, . 
\end{equation}
This definition, while admittedly a little obscure at this point, will
be apparent in the examples presented in Section \ref{sec:examples}.
From \eqref{eq:kappabis} and \eqref{eq:omega}, one obtains:
\begin{equation}
  \label{eq:omegabis}
  \omega_v = 
  \left | \out{v} \right | ,
\end{equation}
and, hence, we can define the bivector $\gw{v}$ as:
\begin{equation}
  \label{eq:Omega}
  \gw{v} = \out{v} \, .
\end{equation}

Based on \eqref{eq:rho} and \eqref{eq:omega} and on the definition of
geometric product given in \eqref{eq:gp}, we can finally provide the
following novel and most important expression of this work:
\begin{equation}
  \label{eq:freq}
  \boxed{\gf{v} = \gr{v}} \, ,
\end{equation}
where we define the term $\gf{v}$ as the \textit{generalized
  frequency} of the voltage $\bfg v$.  The left-hand side of
\eqref{eq:freq} depends only on geometric invariants, namely $v = s'$,
$v' = s''$ and the components of the bivector that define the
magnitude of the curvature $\kappa_v$.  The generalized frequency,
thus, does not depend on the system of coordinates with which $\bfg v$
is represented or measured, nor the number of ``dimensions'' where
$\bfg v$ is defined.  It is also interesting to observe that frequency
is defined in \eqref{eq:freq} as the sum of a symmetric ($\rho_v$) and
an antisymmetric term ($\gw{v}$).  Finally, we note that
\eqref{eq:freq} has been obtained without any assumption on the
dynamic behavior of the components of $\bfg v$.  Unbalanced and/or
non-sinusoidal conditions, multi-phase systems and even dc systems are
consistent with this definition.

\section{Examples}
\label{sec:examples}

The examples presented below are aimed at illustrating the features of
the generalized frequency.  The first two examples show that, in
stationary conditions, \eqref{eq:freq} leads to the well-known and
widely accepted definition of frequency in ac systems.  Examples 3 and
4 illustrate the special cases of transient balanced three-phase
systems and dc systems, respectively.  Example 5 extends the
definition of generalized frequency to the current.  Example
6 illustrates the link between the generalized frequency and the
generalized instantaneous reactive power proposed in \cite{4450060503,
  1308315, 5316097} and shows a simple way to estimate the generalized
frequency in practice.   Finally, Example 7 compares the
estimation of the frequency as obtained with a synchronous reference
frame phase-locked loop (SRF-PLL) and the one obtained with
\eqref{eq:freq} based on a simulation of a detailed EMT model of the
well-known New England 39-bus system.  

\subsubsection*{Example 1}

Let us consider a stationary single-phase voltage with constant angular
frequency $\omega_o$ and magnitude $V$.  The voltage vector can be
defined as:
\begin{equation*}
  \bfg v =
  V \cos(\theta) \e{1} + V \sin(\theta) \e{2} =
  v_1 \e{1} + v_2 \e{2} \, ,
\end{equation*}
where $\theta = \omega_o t + \phi$ and $(\e{1}, \e{2})$ is the
canonical basis of the system, with $\e{1}$ and $\e{2}$
orthonormal vectors.  Then, $v = |\bfg v| = V$ and:
\begin{equation*}
  \bfp v = v'_1\e{1} + v'_2 \e{2} =
  - \omega_o v_2 \e{1} + \omega_o v_1 \e{2} \, ,  
\end{equation*}
from which one can deduce that, as expected:
\begin{align*}
  \rho_v
  &= \frac{1}{V^2}(v_1v'_1 + v_2v'_2) =
    \frac{\omega_o}{V^2}(-v_1v_2 + v_2 v_1) = 0 \, , \\
  \omega_v
  &= \frac{1}{V^2} |v_1v'_2 - v_2v'_1| =
    \frac{\omega_o}{V^2} |v^2_1 + v^2_2| = \omega_o \, .
\end{align*}
It is relevant to note that, in $\mathbb{R}^2$, multivectors are
isomorphic to complex numbers.  In fact, the bivector $\gw{v}$ is:
\begin{equation*}
  \gw{v} = \omega_o (\e{1} \wedge \e{2})  = \jj \omega_o \, .
\end{equation*}
where the imaginary unit $\jj$ is defined as
$\jj \equiv \e{1} \wedge \e{2}$.

\begin{figure}[h!]
  \centering
  \begin{subfigure}{.45\linewidth}
    \centering
    \resizebox{\linewidth}{!}{\includegraphics{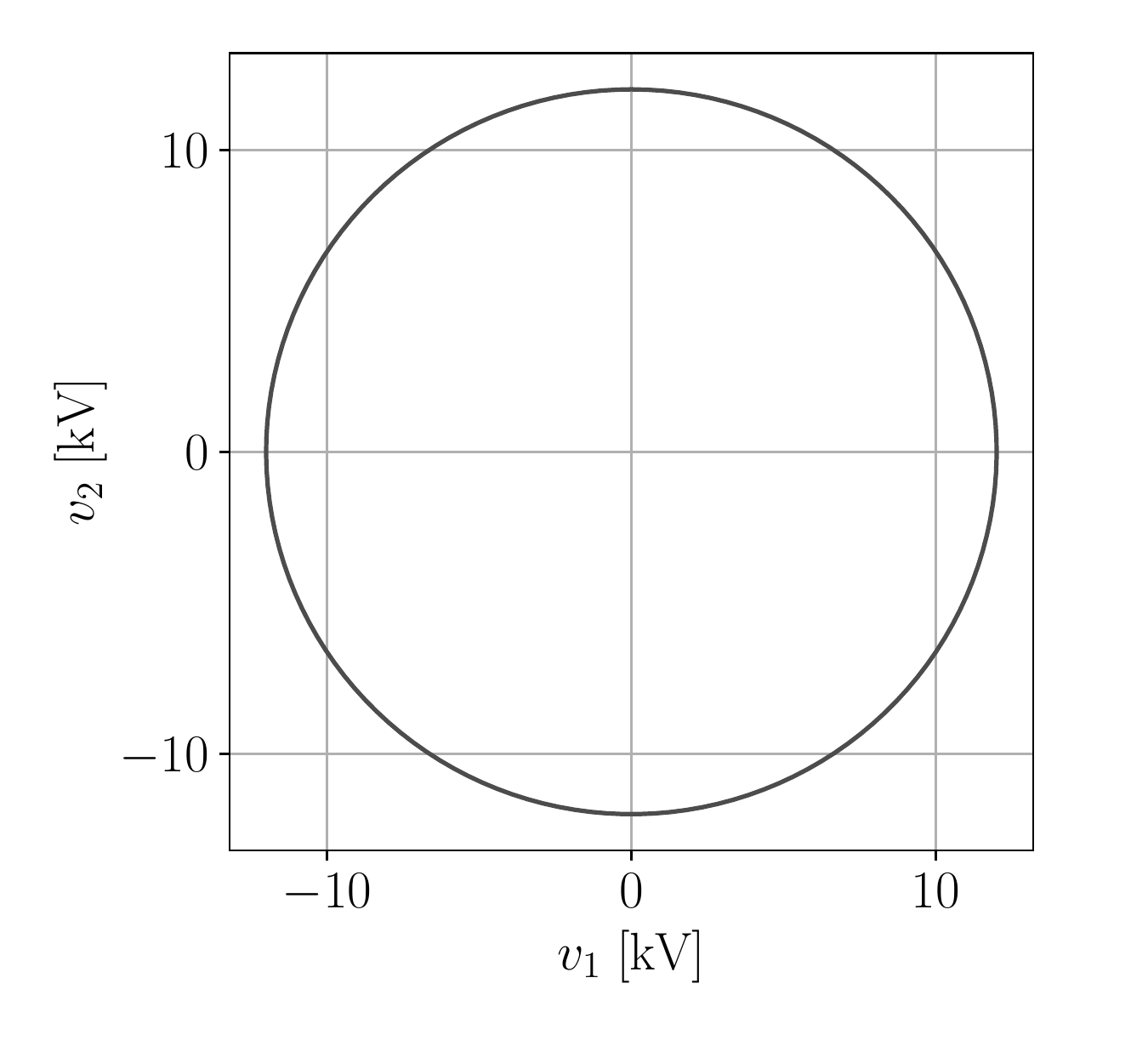}}
    \caption{Example 1}
    \label{fig:E1}
  \end{subfigure}
  \begin{subfigure}{.45\linewidth}
    \centering
    \resizebox{\linewidth}{!}{\includegraphics{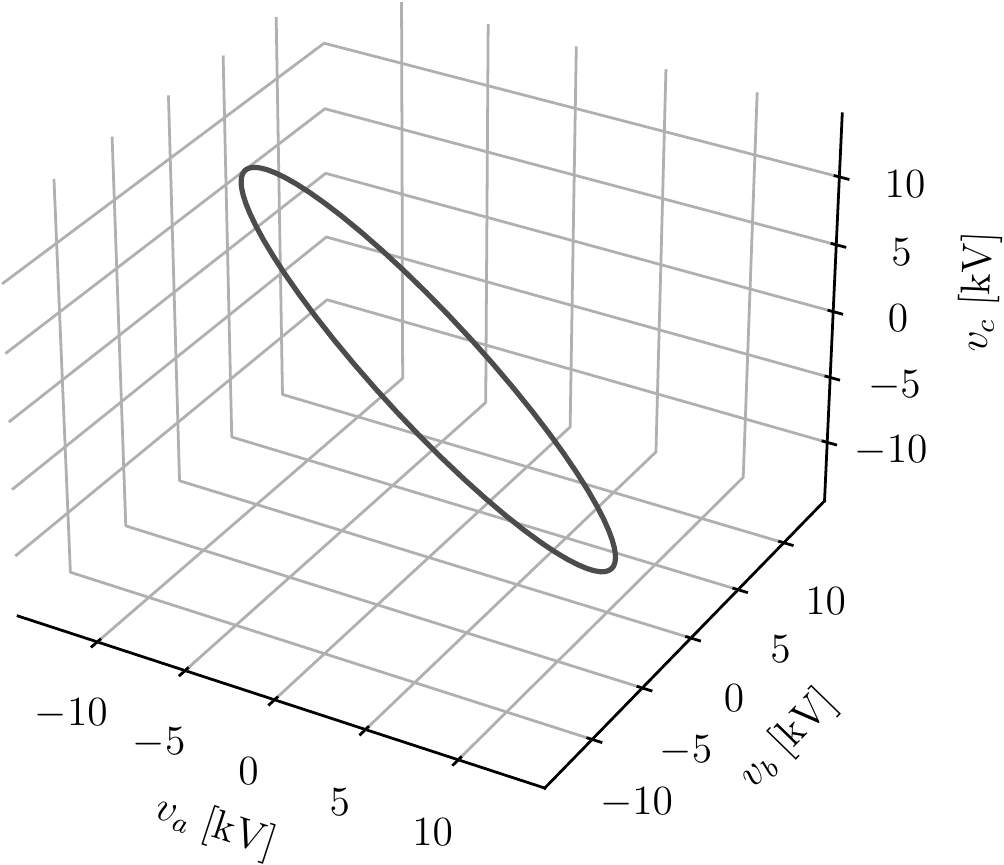}}
    \caption{Example 2}
    \label{fig:E2}
  \end{subfigure}
  \caption{Voltage ``curves.''}
  \label{fig:1}
  \vspace{-2mm}
\end{figure}

It is also worth observing that, in two dimensions, the curvature is
defined as \cite{Stoker}:
\begin{equation}
  \label{eq:theta}
  \kappa_v =
  \frac{d\theta}{ds} =
  \frac{d\theta}{dt}\frac{dt}{ds} =
  \frac{\theta'}{s'} =
  \frac{\theta'}{v} \, ,
\end{equation}
which is valid in any transient and non-sinusoidal conditions.  Using
the chain rule and recalling \eqref{eq:omega}, one has that
$\theta' = \omega$, which is the commonly accepted definition of
frequency \cite{IEEE118}.   Figure \ref{fig:E1}
illustrates the voltage ``curve'' for a single-phase stationary
voltage with $V = 12$ kV and $\omega_o = 120\pi$ rad/s.  As expected,
the curve is a circle, which, as it is well-known, has constant
curvature.

\subsubsection*{Example 2}

We consider a stationary balanced three-phase system.  The voltage
vector is:
\begin{align*}
  \bfg v
  &= V \sin(\theta_a) \e{a} +
    V \sin(\theta_b) \e{b} +
    V \sin(\theta_c) \e{c} \\
  &= v_a \e{a} + v_c \e{b} + v_c \e{c} \, ,
\end{align*}
where $V$ is constant; $\theta_a = \omega_o t$,
$\theta_b = \theta_a - \alpha$ and $\theta_c = \theta_a + \alpha$ with
$\omega_o$ constant and $\alpha = \frac{2\pi}{3}$; and
$(\e{a}, \e{b}, \e{c})$ is the canonical basis of the system, with
$\e{a}$, $\e{b}$ and $\e{c}$ orthonormal vectors.  Then:
\begin{align*}
  v^2 = |\bfg v|^2 = V^2 (\sin^2 \theta_a + \sin^2 \theta_b +
  \sin^2 \theta_c) = \frac{3}{2} V^2 \, ,  
\end{align*}
and
\begin{align*}
  \bfg v'
  &=
    \omega_o V \cos(\theta_a) \e{a} +
    \omega_o V \cos(\theta_b) \e{b} +
    \omega_o V \cos(\theta_c) \e{c} \\
  &= v'_a \e{a} + v'_c \e{b} + v'_c \e{c} \, .
\end{align*}
Then, one has:
\begin{align*}
  \rho_v
  &= \frac{v_a v'_a + v_b v'_b + v_c v'_c}{v^2} \\
  &= \frac{\omega_o V^2 \frac{1}{2}(\sin 2\theta_a +
    \sin 2\theta_b + \sin 2\theta_c) }{v^2} = 0 \, , \\
  \gw{v}
  &= \frac{1}{v^2}
    \begin{bmatrix}
      0 & v_a v'_b - v_b v'_a & v_a v'_c - v_c v'_a \\
      v_b v'_a - v_a v'_b & 0 & v_b v'_c - v_c v'_b \\
      v_c v'_a - v_a v'_c & v_c v'_b - v_b v'_c & 0
    \end{bmatrix} , \\
  \omega_v
  &= \frac{1}{v^2} \sqrt{
    (v_a v'_b - v_b v'_a)^2 +
    (v_b v'_c - v_c v'_b)^2 +
    (v_c v'_a - v_a v'_c)^2} \\
  &= \frac{\sqrt{\omega^2_o V^4(2\sin^2(\alpha) +
    \sin^2(2\alpha))}}{\frac{3}{2} \, V^2} = \omega_o \, .
\end{align*}
Figure \ref{fig:E2} illustrates the voltage ``curve'' for
three-phase balanced and stationary voltages with $V=12$ kV and
$\omega_o = 120\pi$ rad/s.  The curve is a circle in the 3D space and,
as for Example 1, has constant curvature.

\subsubsection*{Example 3}

We consider a balanced three-phase system in transient conditions.  For
illustration, we use the $dqo$ reference frame.  The voltage vector is
$\bfg v = v_d \e{d} + v_q \e{q} + v_o \e{o}$, where
$(\e{d}, \e{q}, \e{o})$ is the canonical basis of the system, with
$\e{d}$, $\e{q}$ and $\e{o}$ orthonormal vectors, and the vectors
$\e{d}$ and $\e{q}$ are rotating at angular speed $\omega_o$.  Since
the system is balanced, $v_o = 0$.  Then, $v^2 = v_d^2 + v_q^2$, and:
\begin{equation*}
  \bfp v = (v'_d - \omega_o v_q) \e{d} +
  (v'_q + \omega_o v_d) \e{q}  =
  \tilde{v}'_d \e{d} + \tilde{v}'_q \e{q} \, ,
\end{equation*}
where, assuming that the $q$-axis leads the $d$-axis, one has:
\begin{equation*}
  \ep{d} = \omega_o \e{q} \, , \quad \ep{q} = - \omega_o \e{d} \, .
\end{equation*}
The components of the generalized frequency are:
\begin{align*}
  \rho_v
  &= \frac{v_d v'_d + \omega_o v_dv_q + v_q v'_q - \omega_o v_dv_q}{v^2} 
    = \frac{v_d v'_d + v_q v'_q}{v^2} \, , \\
  \gw{v}
  &=
    \begin{bmatrix}
      0 & v_d \tilde{v}'_q - v_q \tilde{v}'_d \\
      v_q \tilde{v}'_d - v_d \tilde{v}'_q & 0
    \end{bmatrix} , \\
  \omega_v
  &= \frac{v'_q v_d + \omega_o v^2_d - v'_d v_q + \omega_o v^2_q}{v^2}
  = \omega_o + \frac{v'_q v_d - v'_d v_q}{v^2} \, ,
\end{align*}

\begin{figure}[h!]
  \centering
  \begin{subfigure}{.45\linewidth}
    \centering
    \resizebox{\linewidth}{!}{\includegraphics{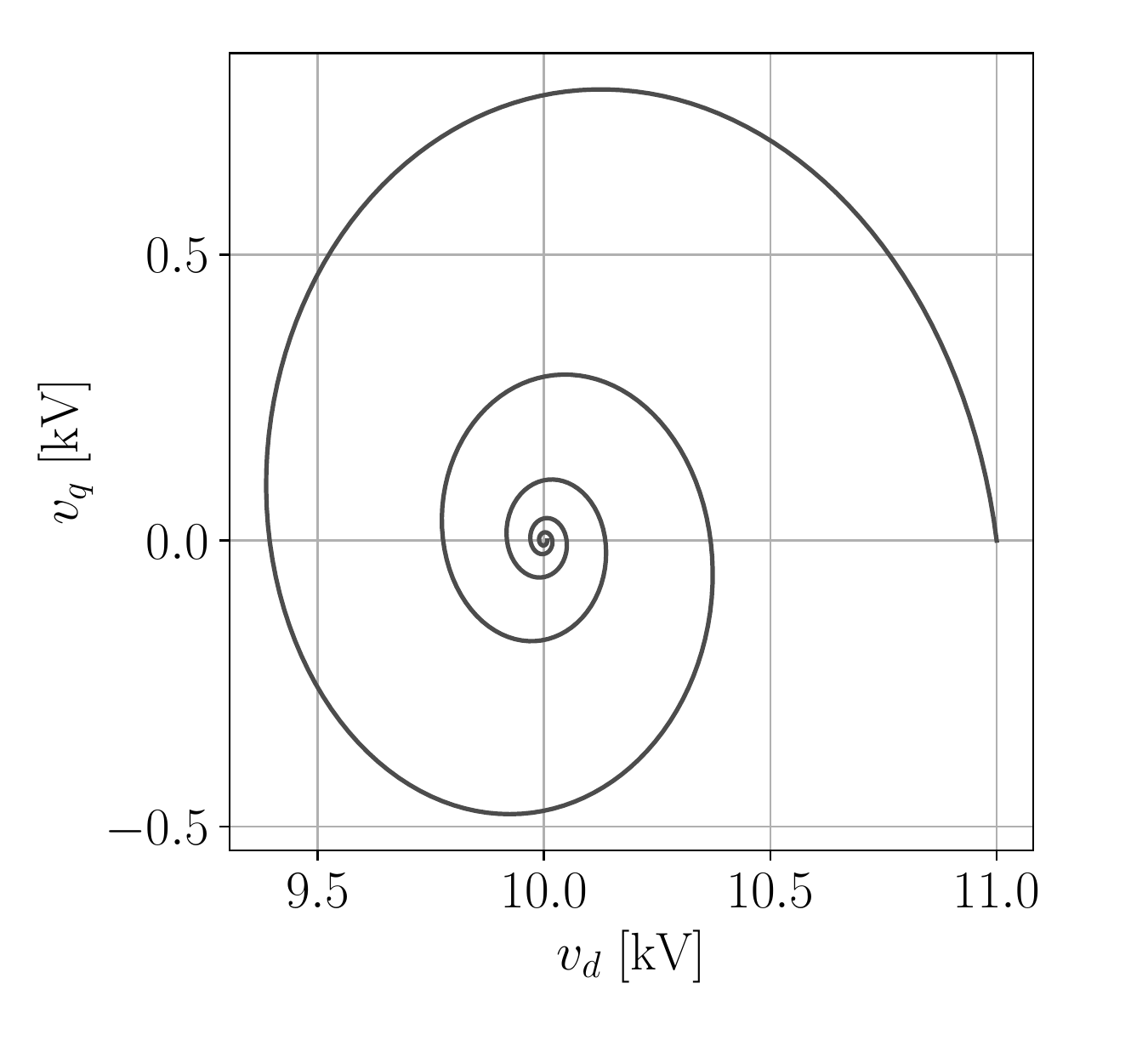}}
    \label{fig:E3.1}
  \end{subfigure}
  \begin{subfigure}{.45\linewidth}
    \centering
    \resizebox{\linewidth}{!}{\includegraphics{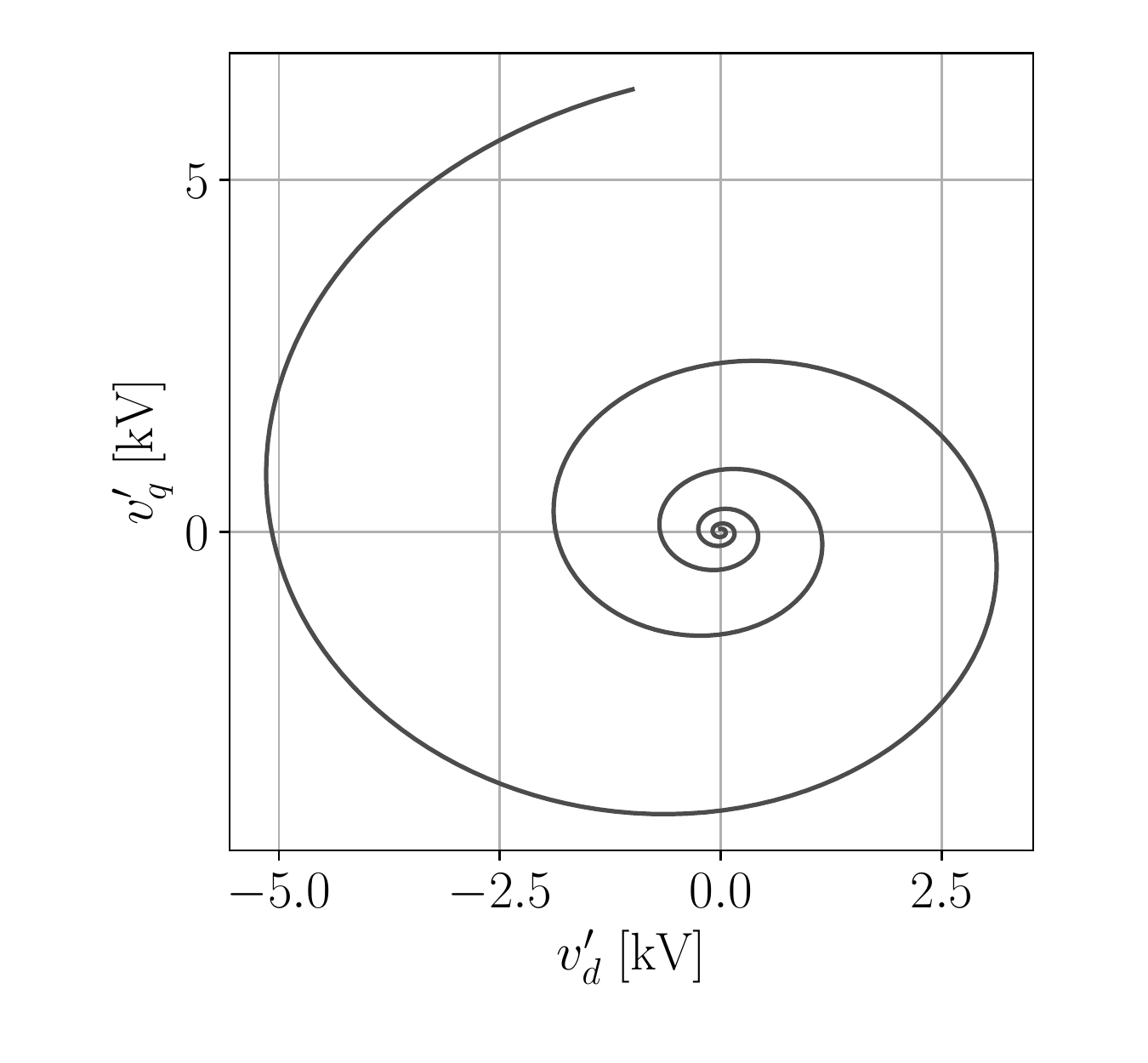}}
    \label{fig:E3.2}
  \end{subfigure}
  \vspace{-5mm} \\
  \begin{subfigure}{.45\linewidth}
    \centering
    \resizebox{\linewidth}{!}{\includegraphics{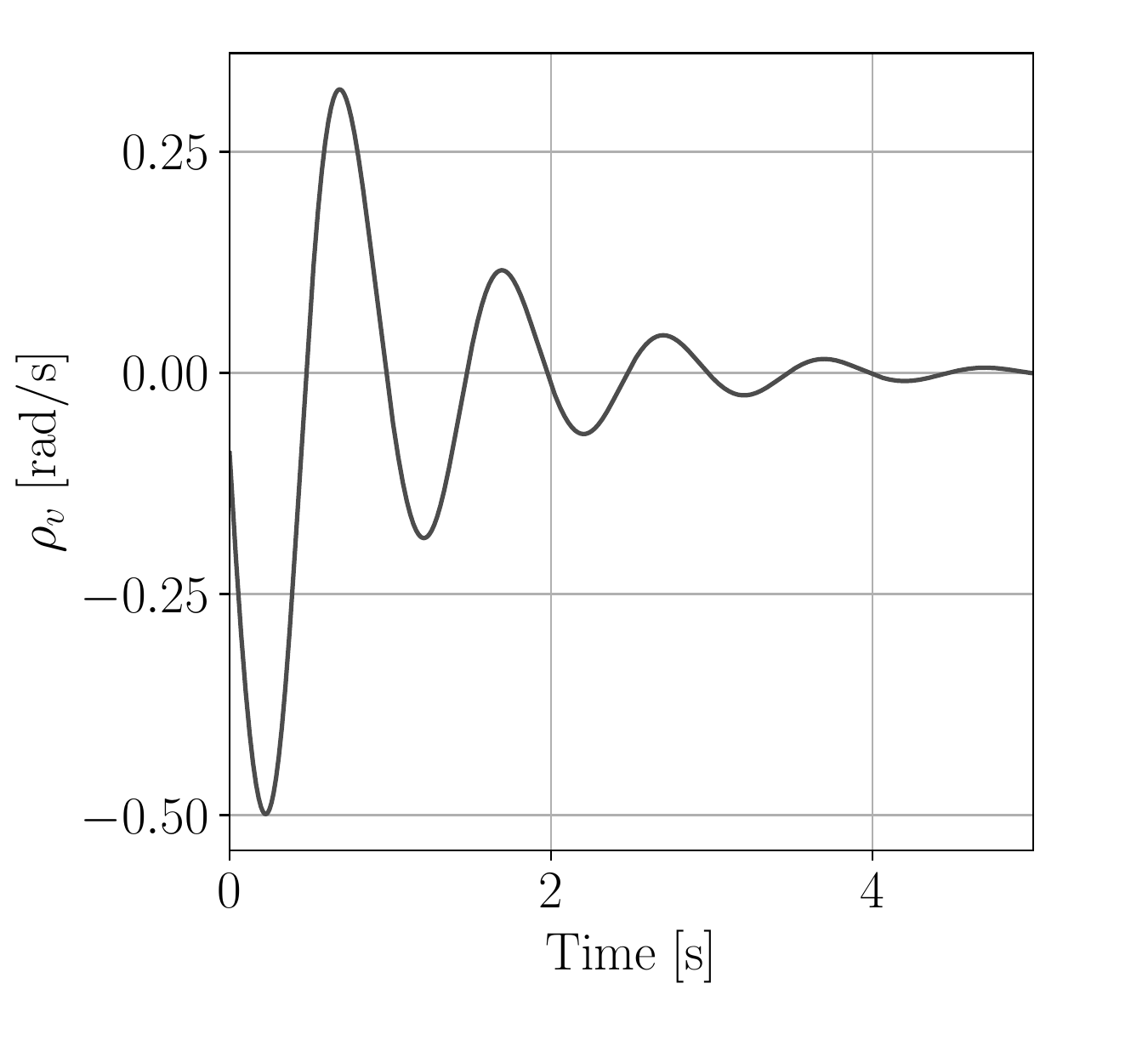}}
    \label{fig:E3.3}
  \end{subfigure}
  \begin{subfigure}{.45\linewidth}
    \centering
    \resizebox{\linewidth}{!}{\includegraphics{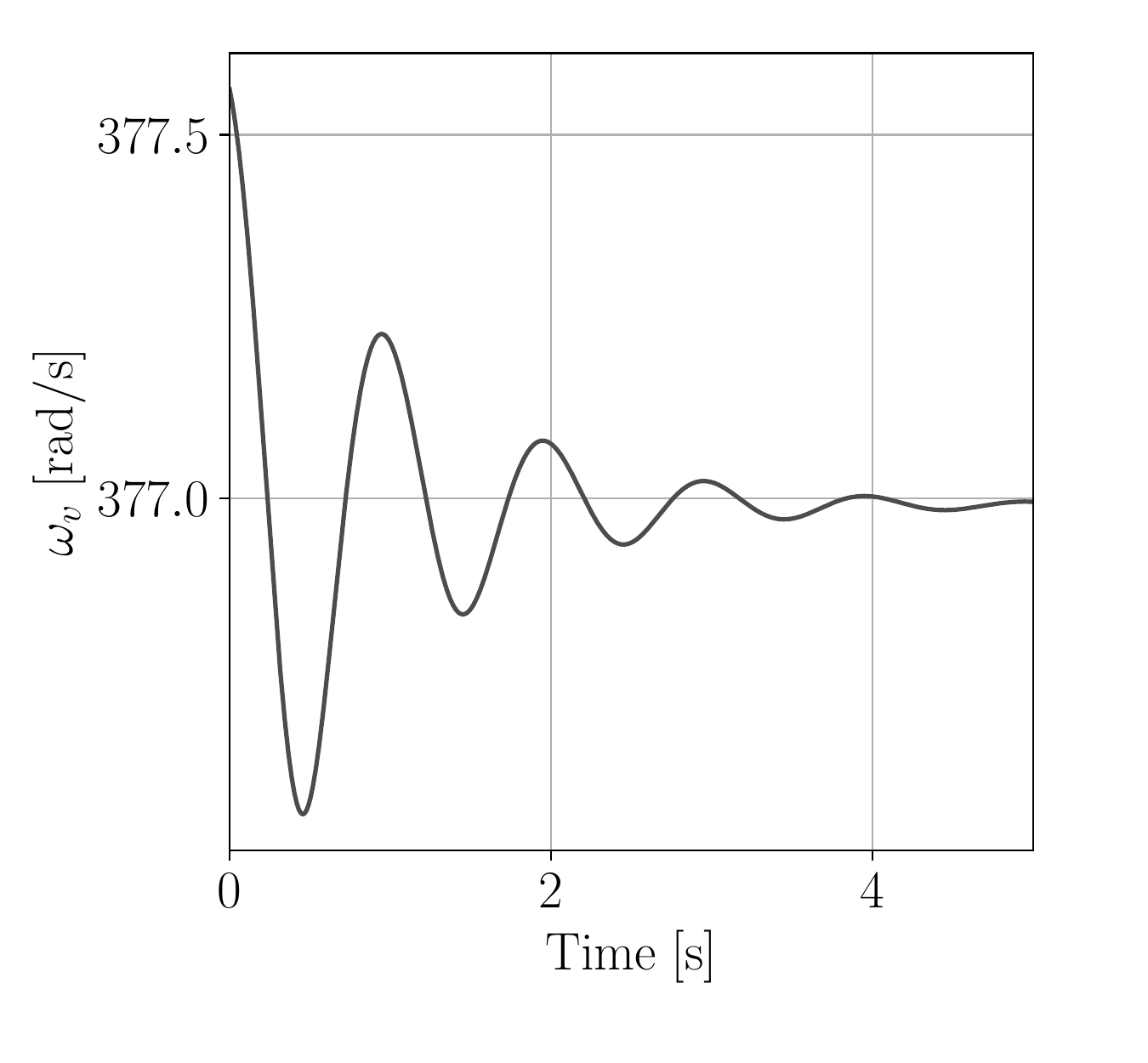}}
    \label{fig:E3.4}
  \end{subfigure}
  \vspace{-5mm}
  \caption{Illustration of Example 3.}
  \label{fig:2}
  \vspace{-2mm}
\end{figure}

The equations above show that the definition of the Park vector as
$\bfg v = v_d + \jj v_q$ and the time derivative in the Park reference
frame, namely $\frac{d}{dt} + \jj \omega_o$, are an equivalent
formulation for balanced three-phase systems in transient conditions
\cite{4450030107}.  Figure \ref{fig:2} illustrates the
expressions above for $\rho_v$ and $\omega_v$ assuming
$v_d = 10 + {\rm exp}(-t)\cos(2\pi t))$ kV,
$v_q = {\rm exp}(-t)\sin(2\pi t))$ kV, $v_o = 0$, and
$\omega_o = 120\pi$ rad/s.  The figure also shows the components of
the voltage and their time derivatives, which confirm that the
curvature and thus the frequency are not constant in this case.

\subsubsection*{Example 4}

We show that \eqref{eq:freq} is valid also for dc voltages.  In dc
circuits, the voltage has only one component along the unique basis of
the system, say $\e{\rm dc}$, hence,
$\bfg v = v_{\rm dc} \, \e{\rm dc}$ and
$\bfp v = v_{\rm dc}' \, \e{\rm dc}$.  From the definitions of inner
and wedge product one has:
\begin{equation*}
  \bfg v \cdot \bfp v = v_{\rm dc} \, v_{\rm dc}' \, , \quad
  \bfg v \wedge \bfp v = \bfg 0 \, .
\end{equation*}
In dc, then, the generalized frequency is equal to
$\rho_v = v_{\rm dc}'/v_{\rm dc}$ and, as expected,
$\omega_v = |\gw{v}| = 0$.

\subsubsection*{Example 5}

Similarly to the voltage, one can define the generalized frequency of
the current.  Consider the vector of the electric charge $\bfg q$ as
an abstract curve in $\mathbb{R}^n$.  This vector does not have to be
intended as a charge moving in space, but rather as the macroscopic
effect of the electric field in a given part of a circuit.  Then:
\begin{equation}
  \label{eq:charge}
  \bfp q = \bfg \ii \, ,
\end{equation}
and, analogously to the discussion on the voltage, the generalized
frequency associated with the current is given by:
\begin{equation}
  \label{eq:current}
  \gf{\ii} = \gr{\ii} \, .
\end{equation}
In general, for any given element of a circuit, the generalized
frequency of the voltage is not equal to that of the current.
Relevant exceptions are resistances.  For a balanced resistive branch,
$\geom{v} = R^2 \geom{\ii}$, which indicates that, from a geometrical
point of view, resistances are scaling factors.

\subsubsection*{Example 6}

We further elaborate on the link between voltage and current vectors.
For balanced capacitive elements, one has:
\begin{equation}
  \label{eq:iCv}
  \bfg \ii = C \bfp v \, ,
\end{equation}
Merging \eqref{eq:freq} and \eqref{eq:iCv}, one obtains:
\begin{equation}
  \label{eq:C}
  \gf{v} =
  \frac{\bfg v \bfg \ii}{Cv^2} =
  \frac{p - \bfb Q}{Cv^2} \, ,
\end{equation}
where $p = \bfg v \cdot \bfg \ii$ is the instantaneous active power,
which is not null only in transient conditions, and
$\bfb Q = \bfg \ii \wedge \bfg v$ is the generalized instantaneous
reactive power as defined in \cite{5316097}.  It is interesting to
note that \eqref{eq:C} provides an expression to calculate the
frequency of an electric circuit in any transient condition through
instantaneous voltage and current measurements.  Interestingly, no
discrete Fourier transforms with mobile windows or other standard
numerical techniques are required.

\subsubsection*{Example 7}

\begin{figure}[h!]
  \centering
  \begin{subfigure}{.49\linewidth}
    \centering
    \resizebox{\linewidth}{!}{\includegraphics{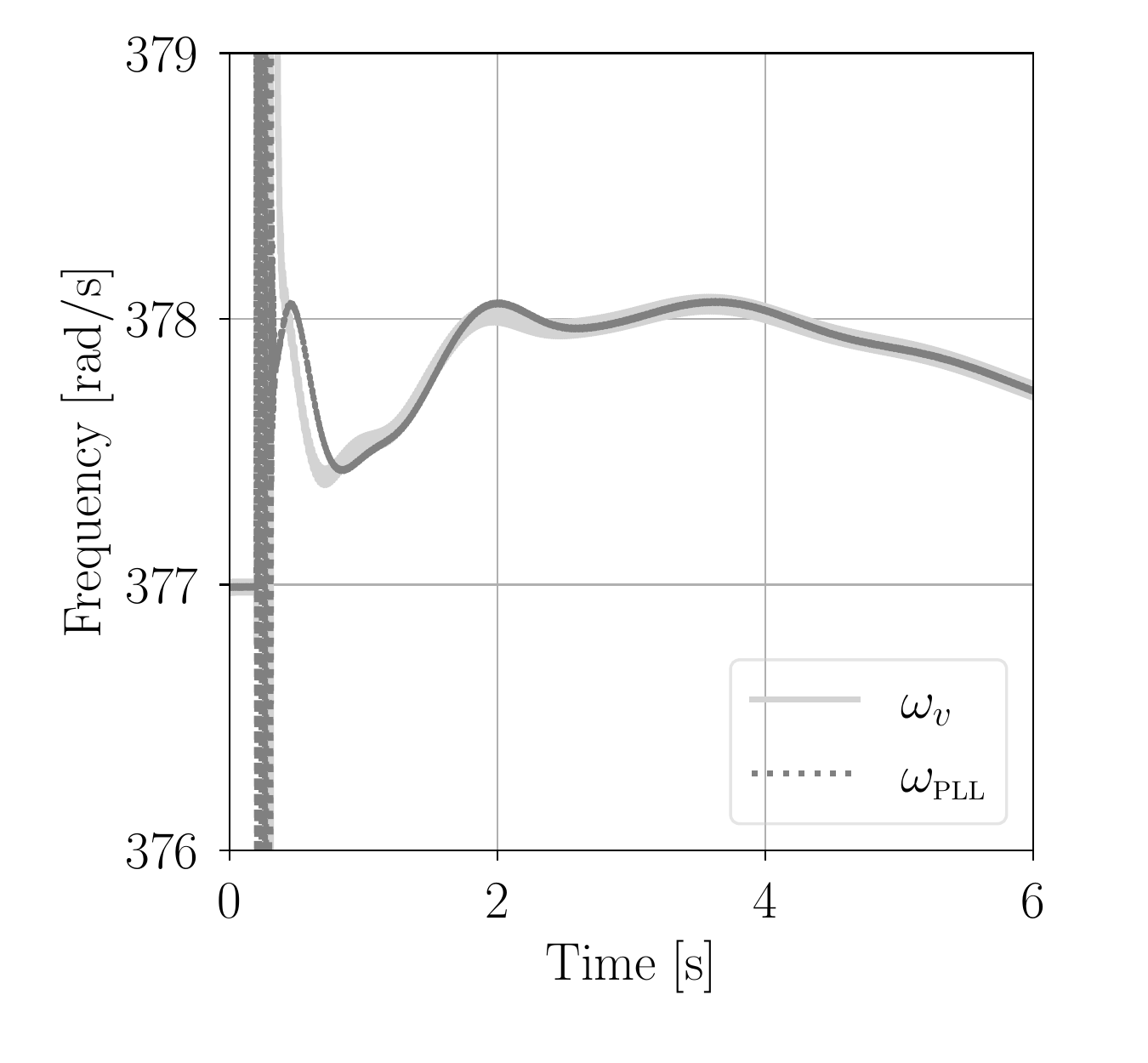}}
  \end{subfigure}
  \begin{subfigure}{.49\linewidth}
    \centering
    \resizebox{\linewidth}{!}{\includegraphics{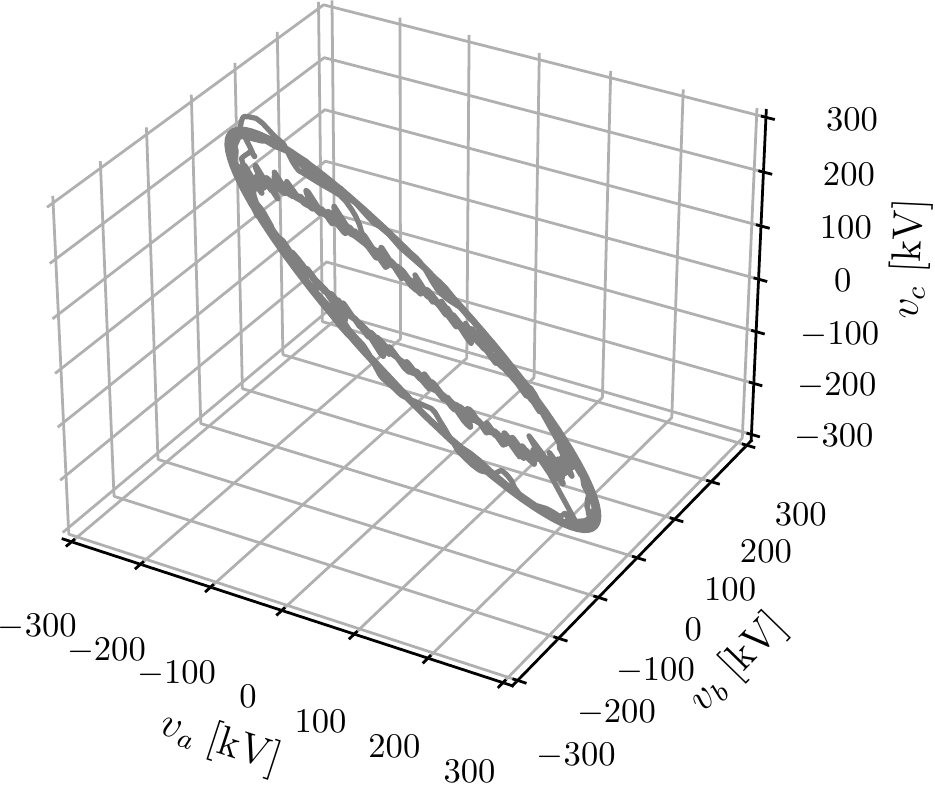}}
  \end{subfigure}
  \caption{Examples 7: Voltage and frequency at bus 26.}
  \label{fig:3}
\end{figure}

This last example presents a comparison of the estimation of the
frequency as obtained with a standard SRF-PLL and with
\eqref{eq:freq}.  Figure \ref{fig:3} shows the results obtained with
PowerFactory and the New England 39-bus system following a
phase-to-phase fault at bus 3 applied at $t = 0.2$ s and cleared at
$t = 0.3$ s.  The simulation utilizes the fully-fledged EMT model
provided by PowerFactory.


PLLs can only estimate $\omega_v = |\bfg \Omega_v|$, i.e., the
magnitude of the bivector defined in \eqref{eq:freq}.  Despite this
limitation of the PLL, the comparison shows that the frequency
obtained with the proposed approach is consistent with the
conventional SRF-PLL.  To obtain the estimation of $\omega_v$ shown in
Fig.~\ref{fig:3} a simple discrete first-order filter is used to
smooth the numerical noise of the time derivative of the voltages and
calculate $\omega_v$.  This is enough in this case to obtain good
results, which are affected by less delay than those obtained with the
PLL.
Finally, the right panel of Fig.~\ref{fig:3} shows that the trajectory
described by the voltage changes plane during the fault and is not
perfectly circular, thus leading to a time-varying curvature/frequency.


\section{Remarks on the Estimation of the Frequency}
\label{sec:remarks}

There exist two main broad approaches for the estimation of the
frequency: transformation-based methods (e.g., Fourier-transform-based
approaches and more recently, Hilbert-transform-based approaches) and
time-domain methods (e.g., PLLs).  The approach proposed in this
letter falls in the second category.
In general, all conventional approaches define \textit{a priori} a
model of the measured signal.  For example, a sine wave with constant
magnitude in \cite{IEEE118}, a set of sine waves with constant
magnitude and frequency in a given window in the case of the Fourier
transform, or a sine wave with pulsating amplitude in the case of the
Hilbert transform \cite{Paolone:2020}.  If the actual signal does not
fit the given model, the estimation obtained with these methods might
not be accurate.

The definition of frequency that is proposed in this letter has two
relevant advantages, as follows.

\begin{itemize}
\item It is intrinsically model agnostic, i.e., no assumption is made
  on the time dependency of the elements of the voltage vector.  This
  allows providing a definition independent from the
  transient/stationary conditions of the circuit where the frequency
  is to be estimated.
\item Equation \eqref{eq:freq} suggests a way to calculate the
  frequency based directly on the measured signal.  Other approaches
  require to process the measurements before being able to do the
  estimation.
\end{itemize}

Note that the issue indicated in the second point above involves also
conventional time-domain approaches, as the signals that feed the PLLs
are the voltage measurements processed through the Park transform,
which is known to be particularly sensitive to noise
\cite{Paolone:2020}.

Of course, also \eqref{eq:freq} presents some challenges.  It requires
in fact to calculate the time derivative of the voltage, which can
lead to numerical issues.  But this is done directly on the measured
quantities not on transformed ones.  Moreover, Example 7 shows that
the issues deriving from numerical differentiation can be resolved
with proper filtering.
  
In summary, the proposed approach appears useful in two ways: (i) to
build a ``theory'' based on the geometrical interpretation of
frequency and electric quantities in general; and (ii) in estimation
and control, to define the frequency in a unambiguous way.  This also
suggests that \eqref{eq:freq}, if standardized, can be utilized to
compare the results obtained with other techniques, e.g., to evaluate
the accuracy of PLLs and other devices that estimate the frequency,
e.g., phasor measurements units.

\section{Conclusions}
\label{sec:conc}


The proposed formal framework generalizes and solves known issues of
the conventional definition of frequency.  A strength of the proposed
approach is that it is based on invariant quantities, hence it is
compatible with any reference frame, e.g., $abc$, $dqo$ and even dc
circuits.  It is interesting to note that the proposed approach
defines the frequency as a geometrical object with symmetric and
antisymmetric parts.  In an example, the letter also shows the link
between the generalized frequency and the power of a circuit.  This
link between geometry and energy appears worth further research.

\section{Acknowledgments}

The author wishes to thank Mr Taulant K{\"e}r{\c c}i and Dr Ahsan
Murad for their help with the simulation results presented in Example
7.



\end{document}